\documentclass[12pt]{amsart}
\usepackage{mathptmx}
\usepackage{amsmath}
\usepackage{amssymb}
\usepackage{array}
\usepackage{geometry}
\usepackage[bookmarks=true,colorlinks=true, pdfstartview=FitV, linkcolor=black, citecolor=blue, urlcolor=black]{hyperref}
\usepackage{color}
\definecolor{DarkRed}{rgb}{0.55,.00,0.2}
\definecolor{DarkGrey}{rgb}{0.35,.35,0.35}

\theoremstyle{definition}

\theoremstyle{remark}

\numberwithin{equation}{section}

\begin{document}

\title{On the iterated  Stieltjes  transform and its convolution with applications to singular integral
equations}

\author{S.Yakubovich and M.Martins}
\address{Department of Mathematics, Fac. Sciences of University of Porto,Rua do Campo Alegre,  687; 4169-007 Porto (Portugal)}
\email{ syakubov@fc.up.pt}



\maketitle

\markboth{\rm \centerline{ S.Yakubovich and M.Martins}}{}
\markright{\rm \centerline{The iterated Stieltjes transform }}

\begin{abstract}
We consider  mapping properties of the iterated  Stieltjes transform, establishing its new relations
with the iterated Hilbert transform (a singular integral) on the half-axis and proving the corresponding convolution and Titchmarsh's type theorems.  Moreover,  the obtained  convolution method is applied to
solve a new class of singular integral equations.
\bigskip

\end{abstract}

\section{Introduction and auxiliary results}

Let $x \in \mathbb{R}_+,\ f  \in  L_p(\mathbb{R}_+),\ 1\le p< \infty$ be a complex-valued function.  It is known that  the classical Laplace transform
$$(\mathcal{ L} f)(x)= \int_0^\infty e^{-xt}f(t)dt,\quad x >0,\eqno(1)$$
is well defined and one can compute its iteration, simply changing the order of integration and calculating an elementary integral. This drives us to the operator of Stieltjes's transform.  Namely, we obtain
$$(\mathcal{S} f)(x) = (\mathcal{ L}^2 f)(x)= \int_0^\infty e^{-xs}\int_0^\infty e^{-st}f(t)dtds
 = \int_0^\infty \frac{f(t)}{x+t} dt,\quad x >0,\eqno(2)$$
where the change of the order of integration is allowed due to Fubini's theorem  via the estimate, which is based on the H$\ddot{o}$lder inequality
$$\int_0^\infty e^{-xs}\int_0^\infty e^{-st} |f(t)| dtds \le \int_0^\infty e^{-xs}\left(\int_0^\infty e^{-q st} dt\right)^{1/q} ds
 \left(\int_0^\infty  |f(t)|^p dt\right)^{1/p}$$
$$= q^{-1/q}  ||f||_p \Gamma \left(1- q^{-1} \right) x^{q^{-1}-1}, \   {1\over q} + {1\over p} =1.$$
Let us compute, in turn, the iteration  of the Stieltjes transform (2) of an arbitrary $f \in  L_p(\mathbb{R}_+),\ 1\le p< \infty$. Similar motivations with the estimate
$$ \int_0^\infty \frac{1}{x+s} \int_0^\infty \frac{|f(t)|}{s+t} dtds \le ||f||_p  \int_0^\infty \frac{1}{x+s}
\left(\int_0^\infty \frac{1}{(s+t)^q} dt\right)^{1/q} ds $$
$$=\left[ \frac{\Gamma(q-1)}{\Gamma(q)}\right]^{1/q} \Gamma(q^{-1})\Gamma(1-q^{-1}) ||f||_p \   x^{q^{-1}-1},
\   {1\over q} + {1\over p} =1 $$
and relations   (2.2.6.24),  (7.3.2.148)  in \cite{prud},  Vol. 1 and  \cite{prud},  Vol.3, respectively,    lead us to the following transformation
$$(\mathcal{S}_2 f)(x) \equiv  G(x)=     \int_0^\infty \frac{\log(x)- \log(t)}{x-t} f(t) dt,\    x >0,\eqno(3)$$
whose  kernel  has a removable singularity at the point $t=x$ and the integral exists in the Lebesgue sense.  This transformation was introduced for the first time in \cite{boas} in the form of the Stieltjes integral and it is called the iterated Stieltjes transform or the $\mathcal{S}_2$-transform.  On spaces of generalized functions the $\mathcal{S}_2$-transform  (3) was extended in \cite{dub} (see also in \cite{prudni}).

In 1990 \cite{con} the first author proposed a new method of convolution constructions for integral transforms, which is based on the double Mellin-Barnes integrals (see in \cite{hai},  \cite{luch}). Following this direction, he established for the first time as an interesting particular case the convolution operator for the Stieltjes transform (2) (see \cite{hai}, formula (24.38))
$$(f *g)_{\mathcal{S}} (x)=  f(x) (Hg)(x)+ g(x) (Hf)(x),\  x >0,\eqno(4)$$
where
$$(Hf)(x)= \int_0^\infty \frac{f(t)}{t- x} dt\eqno(5)$$
is the operator of the Hilbert transform. Moreover, it was proved the corresponding convolution theorem
$$\left( \mathcal{S} \ (f *g)_{\mathcal{S}} \right)(x)= ( \mathcal{S} f)(x)  ( \mathcal{S} g)(x)\eqno(6)$$
in the class of functions, which is associated with the Mellin transform.   Later \cite{sri} these results were extended on
$L_p$-spaces and applied to a class of singular integral equations of convolution type (4).

Our main goal in this paper is to employ  the convolution method to the transformation (3) in order to derive the related
convolution operator, to prove the convolution and Titchmarsh's type theorems (the latter one is about the absence of divisors of zero in the convolution product) and to apply these results, finding solutions and solvability conditions
for a new class of singular integral equations.   However,  we  begin with investigation of mapping properties of
the iterated  Stieltjes  transform (3),  proving  the inversion theorem for this transformation of the Paley-Wiener type,
as it was done by M. Dzrbasjan for the classical Stieltjes transform (2) (see \cite{mhit},  Theorem 2.11) and  which is different from the corresponding inversion in \cite{boas}.

\section{Inversion theorem for the iterated  Stieltjes transformation}

Let $f \in L_2(\mathbb{R}_+)$. Then according to \cite{tit} its
Mellin transform $f^*(s),\ s \in \sigma_s= \{s \in \mathbb{C}: {\rm
Re} s={1\over 2} \}$ is defined by the integral
$$f^*(s)= \int_0^\infty f(t) t^{s-1}dt,\eqno(7)$$
which is convergent in the mean square sense with respect to the
norm in $L_2(\sigma)$. Reciprocally,  the inversion formula takes
place
$$f(x)= {1\over 2\pi i}\int_{\sigma_s} f^*(s) x^{-s} ds\eqno(8)$$
with the convergence of the integral in the mean square sense with
respect to the norm in $L_2(\mathbb{R}_+)$. Furthermore, for any
$f_1,f_2 \in L_2(\mathbb{R}_+)$ the generalized Parseval identity
holds
$$\int_0^\infty f_1(xt)f_2(t) dt= {1\over 2\pi i}\int_{\sigma_s}  f_1^*(s)f_2^*(1-s) x^{-s}
ds, \ x >0\eqno(9)$$
and Parseval's equality of squares of $L_2$- norms
$$\int_0^\infty |f(x)|^2 dx = {1\over 2\pi}
 \int_{-\infty}^{\infty} \left|f^*\left({1\over 2}+ i\tau\right)\right|^2 d\tau.\eqno(10)$$
It is easily seen, that $f(x) \in L_2(\mathbb{R}_+)$ if and only if
${1\over x} f\left({1\over x}\right) \in L_2(\mathbb{R}_+)$.
Hence, writing operator (3) in the form
$$G(x)=     \int_0^\infty \frac{\log(x/t)}{x/t -1} f(t) {dt\over t}=  \int_0^\infty \frac{\log(xt)}{xt -1} f(1/t){dt\over t}$$
and observing that $\log(x)/(x-1) \in L_2(\mathbb{R}_+)$, we appeal
to generalized Parseval equality (9) and relation (8.4.6.11) in
\cite{prud}, Vol. 3 to derive the representation
$$G(x)= \int_0^\infty \frac{\log(xt)}{xt -1} f(1/t){dt\over t}=
{1\over 2\pi i}\int_{\sigma_s}  \left[\Gamma(s)\Gamma(1-s)\right]^2
f^*(s)x^{-s} ds.\eqno(11)$$
As we see from (11), (8) and supplement formula for Euler's
gamma-function, the Mellin transform of $g$ is equal to
$$G^*(s)= \frac{\pi^2}{\sin^2(\pi s)} f^*(s),\quad s \in
\sigma_s.\eqno(12)$$
In order to prove the inversion theorem for transformation (3), we
will employ Dzrbasjan's classes of functions \cite{mhit}, Chap. 2.
Indeed, we have

{\bf Definition 1}. Let $\Phi(s)$ be an entire function
$$\Phi\left(s\right)= \sum_{k=0}^\infty d_k s^k,$$
having on the line $\sigma_s$ the expansion

$$\Phi\left({1\over 2} +i\tau\right)= \sum_{k=0}^\infty
c_k\tau^{2k},\eqno(13)$$
where $c_0 >0, \ c_k \ge 0, \ k=1,2,\dots .$ We will say that $f(x)
\in L_2^{(\Phi)}(\mathbb{R}_+)$, if

 1) $f$ is differentiable infinitely many times on $\mathbb{R}_+$ and
 $$\left(-x{d\over dx}\right)^k f(x) \in L_2(\mathbb{R}_+), \
 k=0,1,2\dots\ ;$$
 2) The following equality holds
 $$\Phi \left(-x{d\over dx}\right)f(x) = \sum_{k=0}^\infty  d_k  \left(-x{d\over dx}\right)^k f(x),$$
where the latter operator series converges in the mean square with
respect to the norm in $L_2(\mathbb{R}_+)$.

Returning to (12), we have
$$f^*(s)= \frac{\sin^2(\pi s)}{\pi^2}G^*(s),\quad s \in \sigma_s.\eqno(14)$$
Hence  letting $\Phi(s)={1\over \pi^2} \sin^2(\pi s)$ in our case (see (13)), we find
$$\Phi\left({1\over 2} +i\tau\right) = \frac{1+\cosh(2\pi\tau)}{2\pi^2} =  {1\over \pi^2} +  {1\over 2\pi^2} \sum_{k=1}^\infty
\frac{(2\pi)^{2k}}{(2k)!} \tau^{2k}.\eqno(15)$$
Now  we are ready to prove the inversion theorem for the iterated  Stieltjes  transform (3).

{\bf Theorem 1}. {\it For an arbitrary function $f \in L_2(\mathbb{R}_+)$ formula $(3)$  defines everywhere on
$\mathbb{R}_+$ a function $G \in  L_2^{(\Phi)}(\mathbb{R}_+)$, where  $\Phi(s)={1\over \pi^2} \sin^2(\pi s)$.
Moreover, almost everywhere the reciprocal inversion formula takes place
$$f(x)=  {1\over \pi^2} \left[ G(x) +  {1\over 2 \sqrt x } \sum_{k= 1}^\infty  \frac{ (-1)^{k} }{(2k)!}
\left(2\pi x{d\over dx}\right)^{2k} \left(\sqrt x \  G(x)\right)\right],\eqno(16)$$
where the  operator series converges in the mean square with respect to the norm in $L_2(\mathbb{R}_+)$.

Conversely, for any $G \in  L_2^{(\Phi)}(\mathbb{R}_+)$ formula $(16)$  defines almost everywhere on
$\mathbb{R}_+$ a function $f \in  L_2(\mathbb{R}_+)$ and the reciprocal formula $(3)$ holds.}

\begin{proof}  According to Definition 1, Parseval equality (10)  and identity  (14), we get   that its right -hand side belongs to $L_2(\sigma_s)$. Thus via Lemma 2.4 from \cite{mhit} $G \in  L_2^{(\Phi)}(\mathbb{R}_+)$. In the meantime evidently for each $k=0,1,2,\dots$
$$\sup_{s \in \sigma_s} \left|s^k \Phi^{-1}(s)\right| = a_k < \infty.$$
Hence due to inequalities
$$\int_{\sigma_s} \left|s^k G^*(s)\right| ^2 |ds| \le a^2_k \int_{\sigma_s} \left|\Phi(s) G^*(s)\right| ^2 |ds|,
\  k=0,1,2,\dots $$
we have $s^k  G^*(s) \in L_2(\sigma_s) \cap L_1(\sigma_s), \   k=0,1,\dots\ .$ Hence owing to the differentiation under the integral sign we find immediately the  representations
$$\left(x{d\over dx}\right)^{2k} \left(\sqrt x \  G(x)\right) =
{1\over 2\pi i} \int_{\sigma_s} (s- 1/2) ^{2k} G^*(s) \  x^{1/2 -s} ds,\ k=0,1,\dots,$$
and therefore,
 $${1\over \sqrt x} P_n \left(x{d\over dx}\right) \left(\sqrt x\ G(x)\right) = {1\over 2\pi i} \int_{\sigma_s} P_n(s) G^*(s) x^{-s} ds,\eqno(17)$$
where by $P_n(s)$ we denoted the sum
$$P_n(s) =  {1\over \pi^2}+  {1\over 2\pi^2} \sum_{k= 1}^n
(-1)^{k} \frac{(2\pi)^{2k}}{(2k)!} (s-1/2)^{2k}.$$
But (14) and (8) yield
$$f(x)= {1\over 2\pi i} \int_{\sigma_s} \Phi(s) G^*(s) x^{-s} ds$$
and from (17) it has for each $n \in \mathbb{N}$
$$f(x)- {1\over \sqrt x} P_n \left(x{d\over dx}\right) \left(\sqrt x \  G(x)\right)=  {1\over 2\pi i} \int_{\sigma_s} \Phi(s) \left[1- \frac{P_n(s)}{\Phi(s)}\right] G^*(s) x^{-s} ds.$$
Therefore  appealing to the Parseval equality (10) we derive
$$\int_0^\infty \left| f(x)- {1\over \sqrt x} P_n \left(x{d\over dx}\right) \left(\sqrt x \  G(x)\right) \right|^2 dx= {1\over 2\pi} \int_{\sigma_s}
 |\Phi(s) G^*(s)|^2 \left|1- \frac{P_n(s)}{\Phi(s)}\right|^2 | ds|.$$
However, the right-hand side of the latter equality tends to zero, when $n \to \infty$ by virtue of the Lebesgue dominated convergence theorem,   since
$$\lim_{n\to \infty} \left|1- \frac{P_n(s)}{\Phi(s)}\right|= 0, \quad s \in \sigma_s$$
and via (15)  for all $n$
$$ \left|1- \frac{P_n(s)}{\Phi(s)}\right| \le 2,   \quad s \in \sigma_s.$$
Thus we arrive at the inversion formula (16), where the operator series converges in the mean square
with respect to the norm in $L_2(\mathbb{R}_+)$.   In the same way, starting from (16) and using (12),  we prove the converse proposition of the theorem.

\end{proof}

\section{The convolution and Titchmarsh theorems}

In this section we will construct and study mapping properties of
the convolution, related to the iterated  Stieltjes   transform
(3). In fact, according to  formula (12.22) in \cite{luch}  of the
generalized $G$-convolution we have

{\bf Definition 2}. Let $f, g$ be functions from $\mathbb{R}_+$ into $\mathbb{C}$ and $f^*,\  g^*$ be their Mellin transforms $(7)$.  Then the function $f*g$ being defined on $\mathbb{R}_+$ by the double Mellin-Barnes integral
$$(f*g)(x)= \frac{\sqrt x}{(2\pi i)^2} \int_{\sigma_s}  \int_{\sigma_w}\left[
 \frac{ \Gamma(s)\Gamma(1-s) \Gamma(w)\Gamma(1-w)}
 { \Gamma(s+w-1/2)\Gamma(3/2-s-w)}\right]^2 $$$$\times   f^*(s)g^*(w) x^{-s-w} dsdw\eqno(18)$$
  is called the convolution of $f$ and $g$ (provided that it exists).

  Using again the supplement formula for gamma-functions and elementary trigonometric identities,  we obtain
  $$ \frac{ \Gamma(s)\Gamma(1-s) \Gamma(w)\Gamma(1-w)}
 { \Gamma(s+w-1/2)\Gamma(3/2-s-w)}=  \pi \left[1- \cot(\pi s)\cot(\pi w)\right].$$

{\bf Lemma 1}. {\it Let $f, g$ be such that their Mellin transforms $f^*,\ g^*$ satisfy conditions $s f^*(s),\  s g^*(s)  \in  L_2(\sigma_s)$.
 Then convolution $(18)$ exists as a continuous function on $\mathbb{R}_+$,  $f*g \in L_2(\mathbb{R}_+)$  and the following inequalities  hold}
$$|(f*g)(x)| \le {2\pi \over \sqrt x}  \left(\int_{-\infty}^\infty  \left|\left({1\over 2}+ i\tau\right) f^*\left({1\over 2}+ i\tau\right)\right|^2
d\tau\right)^{1/2}$$$$\times \left(\int_{-\infty}^\infty  \left|\left({1\over 2}+ i\theta\right) g^*\left({1\over 2}+ i\theta\right)\right|^2 d\theta
\right)^{1/2},\eqno(19)$$
$$\int_0^\infty |(f*g)(x)|^2 dx  \le 16\pi^2  \int_{-\infty}^\infty  \left|  \left({1\over 2}+ i\theta \right)
 f^*\left({1\over 2}+ i\theta\right) \right|^2 d\theta $$$$\times  \int_{-\infty}^\infty
  \left| \left({1\over 2}+ i\tau \right) g^*\left({1\over 2}+ i\tau\right)\right|^2 d\tau .\eqno(20)$$
\begin{proof} In fact,  with the Schwarz inequality for double integrals,   inequality $|\tanh(\pi\tau)| \le 1, \  \tau \in \mathbb{R}$ and
computation of elementary integrals, we obtain
$$|(f*g)(x)| \le {1\over 4 \sqrt x}  \int_{-\infty}^\infty  \int_{-\infty}^\infty
 \left[\tanh(\pi\tau)\tanh(\pi\theta) +1\right]^2 $$$$\times
  \left| f^*\left({1\over 2}+ i\tau\right)g^*\left({1\over 2}+ i\theta\right) \right|d\tau d\theta\le
  {1\over 4 \sqrt x}  \left(\int_{-\infty}^\infty  \int_{-\infty}^\infty \frac{
 \left[\tanh(\pi\tau)\tanh(\pi\theta) +1\right]^2 }{\theta^2+1/4} \right.$$$$\left.\times
  \left|\left({1\over 2}+ i\tau\right) f^*\left({1\over 2}+ i\tau\right)\right|^2d\theta d\tau\right)^{1/2}
  \left(\int_{-\infty}^\infty  \int_{-\infty}^\infty \frac{ \left[\tanh(\pi\tau)\tanh(\pi\theta) +1\right]^2 }{\tau^2+1/4}
   \right.$$$$\times \left. \left|\left({1\over 2}+ i\theta\right) g^*\left({1\over 2}+ i\theta\right)\right|^2d\theta d\tau\right)^{1/2}
   \le {2\pi\over \sqrt x}  \left(\int_{-\infty}^\infty  \left|\left({1\over 2}+ i\tau\right)
    f^*\left({1\over 2}+ i\tau\right)\right|^2 d\tau\right)^{1/2}
    $$$$\times \left(\int_{-\infty}^\infty  \left|\left({1\over 2}+ i\theta\right) g^*\left({1\over 2}+ i\theta\right)\right|^2 d\theta\right)^{1/2},$$
which leads to (19) and guarantees continuity of the convolution
$(f*g)(x)$ on $\mathbb{R}_+$ via the Weierstrass test of the uniform
convergence of the double integral (18) for $x \ge x_0 >0$.
Furthermore, appealing to the Parseval equality (10) and making a
simple change of variables $z= s+w-1/2$ in (18), we get
$$\int_0^\infty |(f*g)(x)|^2 dx = {\pi\over 8}  \int_{-\infty}^\infty
\left| \int_{-\infty}^\infty \left[\tanh(\pi\theta)\tanh(\pi(\tau-\theta)) +1\right]^2 \right.
$$$$\left. \times f^*\left({1\over 2}+ i\theta\right)g^*\left({1\over 2}+ i(\tau- \theta)\right) d\theta \right|^2 d\tau
\le 2 \pi \int_{-\infty}^\infty  d\tau $$$$\times
\left|\int_{-\infty}^\infty   f^*\left({1\over 2}+
i\theta\right)g^*\left({1\over 2}+ i(\tau- \theta)\right) d\theta
\right|^2.$$
Hence   we employ the generalized Minkowski inequality to derive
 $$\left(\int_0^\infty |(f*g)(x)|^2 dx\right)^{1/2} \le  \sqrt{2\pi} \int_{-\infty}^\infty  \left| f^*\left({1\over 2}+ i\theta\right) \right|
  \left( \int_{-\infty}^\infty  \left| g^*\left({1\over 2}+ i(\tau- \theta)\right)\right|^2 d\tau \right)^{1/2} d\theta$$
$$\le 2\sqrt{2 \pi} \int_{-\infty}^\infty  \left| f^*\left({1\over 2}+ i\theta\right) \right| d\theta
\left( \int_{-\infty}^\infty  \left| \left({1\over 2}+ i\tau \right)
g^*\left({1\over 2}+ i\tau\right)\right|^2 d\tau \right)^{1/2} $$
 $$\le  2\sqrt{2 \pi}  \left( \int_{-\infty}^\infty  {d\theta\over \theta^2+ 1/4}\right)^{1/2}
\left( \int_{-\infty}^\infty  \left|  \left({1\over 2}+ i\theta
\right) f^*\left({1\over 2}+ i\theta\right) \right|^2
d\theta\right)^{1/2}
$$$$\times  \left( \int_{-\infty}^\infty  \left| \left({1\over 2}+ i\tau \right) g^*\left({1\over 2}+ i\tau\right)\right|^2 d\tau \right)^{1/2} $$
 $$=  4 \pi \left( \int_{-\infty}^\infty  \left|  \left({1\over 2}+ i\theta \right)
  f^*\left({1\over 2}+ i\theta\right) \right|^2 d\theta\right)^{1/2}
  \left( \int_{-\infty}^\infty  \left| \left({1\over 2}+ i\tau \right)
   g^*\left({1\over 2}+ i\tau\right)\right|^2 d\tau \right)^{1/2} .$$
Thus squaring both sides of the inequalities,  we arrived at (20) and proved the lemma.
 \end{proof}

 Now we are ready to prove the convolution theorem for transformation (3).  Precisely, we state

 {\bf Theorem 2}.  {\it  Let $f^*,\ g^*$ be the Mellin transforms of $f, g$, respectively, satisfying
  conditions $s f^*(s),\ s g^*(s)  \in  L_2(\sigma_s)$.   Then the Mellin transform of the convolution
  $(18)\    (f*g)^*(s) \in L_2(\sigma_s)$ and is equal to
 $$(f*g)^*(s) =   \frac{\pi}{2i}  \int_{\sigma_w}    \left[1+ \tan(\pi( s-w))\cot(\pi w)\right]^2  f^*(s-w+1/2 )g^*(w)dw.\eqno(21)$$
 Moreover,  the factorization equality  holds
$$(\mathcal{S}_2 (f*g) )(x) = \sqrt x \ (\mathcal{S}_2 f)(x) (\mathcal{S}_2 g)(x), \quad x >0. \eqno(22)$$
Besides, if  $s f^*(s), \  s g^*(s)   \in  L_2(\sigma_s) \cap  L_1(\sigma_s)$, then for all $x >0 $ the following representation takes place
$$(f*g)(x)= \pi^2 \sqrt x \left[  f(x) g(x)   - {2\over \pi^2} (Hf)(x) (Hg)(x) \right.$$$$\left. + {1\over \pi^4} (H^2 f)(x) (H^2 g)(x)\right],\eqno(23) $$
where $H$ is the operator of the Hilbert transform $(5)$ and by $H^2$ the iterated Hilbert transform is denoted.}

\begin{proof}   Formula (21) and condition $(f*g)^*(s) \in L_2(\sigma_s)$ follow immediately from (18),
Lemma 1 and Parseval's  equality (10).  Hence via generalized Parseval's identity (9) (see also (3), (11)) it has
$$(\mathcal{S}_2 (f*g) )(x) =   \frac{1}{(2\pi i)^2} \int_{\sigma_s} {\pi^4 x^{-s} \over \sin^2(\pi s) }
 \int_{\sigma_w}    \left[1+ \tan(\pi( s-w))\cot(\pi w)\right]^2$$$$
 \times   f^*(s-w+1/2 )g^*(w)\  dw ds= \frac{1}{(2\pi i)^2} \int_{\sigma_s} x^{-s}
 \int_{\sigma_w}      \frac{ \pi^4 f^*(s-w+1/2 )g^*(w)} {\sin^2(\pi(s-w+1/2))\sin^2(\pi w)} \   dw ds$$
$$=  \frac{\sqrt x }{(2\pi i)^2} \int_{\sigma_s}
 \frac{ \pi^2 f^*(s )} {\sin^2(\pi s)}  x^{-s}\  ds \int_{\sigma_w}
   \frac{\pi^2 g^*(w)} {\sin^2(\pi w)}  x^{-w}\   dw =  \sqrt x \ (\mathcal{S}_2 f)(x) (\mathcal{S}_2 g)(x), \quad x >0,$$
where the change of the order of integration is by Fubini's theorem by virtue of the estimate
$$\int_{\sigma_s}  \int_{\sigma_w}    \left|  \frac{ f^*(s-w+1/2 )g^*(w)} {\sin^2(\pi(s-w+1/2))\sin^2(\pi w)} \right| \ | dw ds|$$$$
= \int_{-\infty}^\infty   \int_{-\infty}^\infty     \left|  \frac{ f^*(i(\tau-\theta) +1/2 )g^*(i \theta+1/2)} {\cosh^2(\pi(\tau-\theta))
\cosh^2(\pi \theta)} \right| \  d\tau  d\theta \le 4 \left( \int_{-\infty}^\infty
   \left| ( i\tau +1/2 ) f^*(i\tau +1/2 )\right|^2 d\tau\right)^{1/2}
   $$$$\times  \left( \int_{-\infty}^\infty  \left| ( i\theta  +1/2 )  g^*(i \theta+1/2)\right|^2 d\theta \right)^{1/2} < \infty. $$
So we established equality (22). In order to prove (23), we return to (18) and write it in the form
$$(f*g)(x)= \frac{\sqrt x}{(2\pi i)^2} \int_{\sigma_s} \int_{\sigma_w}\left[1-  2\cot(\pi s)\cot(\pi w) + \cot^2(\pi s)\cot^2(\pi w)
  \right] $$$$\times   \pi^2 f^*(s)g^*(w) x^{-s-w} dsdw =  \pi^2 \sqrt x \  \left[  f(x) g(x)   +   \frac{1}{2\pi^2}
   \int_{\sigma_s}  \cot(\pi s)  f^*(s) x^{-s} ds \right.$$$$\left.
   \times \int_{\sigma_w} \cot(\pi w) \ g^*(w) x^{-w} dw -  \frac{1}{4\pi^2} \int_{\sigma_s}  \cot^2(\pi s)
     f^*(s) x^{-s} ds \right.$$$$\left.\times \int_{\sigma_w} \cot^2(\pi w) \ g^*(w) x^{-w} dw\right]\eqno(24)$$
and the latter equality in (24) is indeed possible since owing to conditions of the theorem $f^*(s),\ g^*(s)  \in  L_1(\sigma_s)$.
 Now our goal is to prove the equalities
$$\frac{1}{2\pi i} \int_{\sigma_s}  \cot(\pi s)  f^*(s) x^{-s} ds= {1\over \pi} (Hf)(x),\quad x >0,\eqno(25)$$
$$\frac{1}{2\pi i} \int_{\sigma_s}  \cot^2(\pi s)  f^*(s) x^{-s} ds= {1\over \pi^2} (H^2 f)(x),\quad x >0.\eqno(26)$$
In order to do this, we employ the known equality (\cite{prud},  Vol. 1, relation (2.2.4.26))
$$PV \  {1\over \pi} \int_0^\infty {t^{s-1}\over 1-t}\ dt = \cot(\pi s), \quad   0 < {\rm Re} s < 1,\eqno(27)$$
where its   left-hand side  is understood as
$$PV \   {1\over \pi} \int_0^\infty {t^{s-1}\over 1-t}\ dt= \lim_{\varepsilon \to 0} \varphi_\varepsilon(s),$$
and
$$\pi \varphi_\varepsilon (s)=  \left(  \int_0^{1-\varepsilon} +  \int_{1+\varepsilon}^\infty  \right) {t^{s-1}\over 1-t}\ dt,  \quad 0 < \varepsilon < 1, \   0 < {\rm Re} s < 1.\eqno(28)$$
We will treat the following integral
$$I_\varepsilon(x)= \frac{1}{2\pi i} \int_{\sigma_s}  \varphi_\varepsilon (s)  f^*(s) x^{-s} ds,\eqno(29)$$
showing,  that it is possible to pass to the limit under the integral sign when $\varepsilon \to 0$.  This fact will be done, establishing the uniform estimate
$$ \left|\varphi_\varepsilon(s)\right| \le C |s|,\quad  s \in \sigma_s,\eqno(30)$$
where $C >0$ is an absolute constant.  So, in order to prove (29), we choose $0 < \varepsilon < 1/2$ and split integrals in (28) as follows
$$\pi \varphi_\varepsilon (s)=  \left( \int_0^{1/2}+   \int_{1/2}^{1-\varepsilon} + \int_{1+\varepsilon}^{3/2} +
 \int_{3/2}^\infty  \right) {t^{s-1}\over 1-t}\ dt$$$$=  I_1(s)+ I_2(s)+ I_3(s)+ I_4(s),  \quad s \in \sigma_s.$$
Clearly,
$$ \left|I_1(s)\right| \le  \int_0^{1/2}{dt\over (1-t)\sqrt t }= O(1).$$
 Analogously,
$$ \left|I_4(s)\right| \le  \int_{3/2}^\infty {dt\over (t-1)\sqrt t }= O(1).$$
Concerning integral $I_2$, we have $(s= 1/2 +i\tau)$
$$I_2(s)=  \int_{1/2}^{1-\varepsilon}  {\cos(\tau\log t) + i\sin(\tau\log t)\over (1-t)\sqrt t}\ dt $$
and via elementary inequality $|\sin x | \le |x|, \ x \in \mathbb{R}$
$$\left| \int_{1/2}^{1-\varepsilon}  {\sin(\tau\log t )\over (1-t)\sqrt t}\ dt \right| \le |\tau| \int_{1/2}^{1}
 {|\log t | \over (1-t)\sqrt t}\ dt =O(\tau).$$
Further,
$$\int_{1/2}^{1-\varepsilon}  {\cos(\tau\log t)\over (1-t)\sqrt t}\ dt = \int_{1/2}^{1-\varepsilon}  {\cos(\tau\log t) -1\over (1-t)\sqrt t}\ dt
 + \int_{1/2}^{1-\varepsilon}  {1\over (1-t)\sqrt t}\ dt $$
and after integration by parts in the second integral, we find
$$ \int_{1/2}^{1-\varepsilon}  {1\over (1-t)\sqrt t}\ dt = -{ \log \varepsilon\over \sqrt{1-\varepsilon}}-  \sqrt 2 \log 2-
 {1\over 2} \int_\varepsilon^{1/2}  {\log t\over (1-t)^{3/2}}\ dt.$$
In the meantime, with the Lagrange theorem
$$ {\cos(\tau\log t) -1\over t-1} =   -  \tau \   {\sin(\tau\log (\xi_t)) \over \xi_t},\quad 1/2\le   t< \xi_t< 1. $$
Hence,
$$\left| \int_{1/2}^{1-\varepsilon}  {\cos(\tau\log t) -1\over (1-t)\sqrt t}\ dt \right| \le 2 |\tau| \int_{1/2}^{1}  {dt \over \sqrt t}
= O(\tau).$$
Similarly,
$$I_3(s)=  \int_{1+\varepsilon}^{3/2}  {\cos(\tau\log t) + i\sin(\tau\log t)\over (1-t)\sqrt t}\ dt $$
and
$$\left| \int_{1+\varepsilon}^{3/2}  {\sin(\tau\log t)\over (1-t)\sqrt t}\ dt\right|\le  |\tau| \int_{1}^{3/2}
 {|\log t | \over (t-1)\sqrt t}\ dt =O(\tau).$$
Meanwhile,
$$\int_{1+\varepsilon}^{3/2}   {\cos(\tau\log t)\over (1-t)\sqrt t}\ dt
 = \int_{1+\varepsilon}^{3/2}   {\cos(\tau\log t)-1 \over (1-t)\sqrt t}\ dt  + \int_{1+\varepsilon}^{3/2}   {1\over (1-t)\sqrt t}\ dt  $$
and , in turn, with the same arguments
$$ \int_{1+\varepsilon}^{3/2}   {1\over (1-t)\sqrt t}\ dt  =  { \log \varepsilon\over \sqrt{1+\varepsilon}} + \sqrt {2/3}  \log 2-
 {1\over 2} \int_{\varepsilon}^{1/2}  {\log t\over (1+t)^{3/2}}\ dt,$$
$$\left| \int_{1+\varepsilon}^{3/2}   {\cos(\tau\log t)-1 \over (1-t)\sqrt t}\ dt \right| \le  |\tau| \int_{1}^{3/2}  {dt \over \sqrt t}
= O(\tau).$$
Thus,
$$\left| I_2(s)+I_3(s)\right| \le 2 \varepsilon |\log\varepsilon| + O(1)+ O(\tau)<  \log 2 +  O(1)+ O(\tau)$$
and combining with estimates above,  we complete the proof of  inequality (30). Returning to (29) and appealing to the Lebesgue dominated convergence theorem, one can pass to the limit when $\varepsilon \to 0$ under the integral sign.
Consequently,  employing (28) and making  simple changes of variables by Fubini's theorem $(f^* \in L_1(\sigma_s))$ with  the use of (8), we obtain  for all $x>0$
$$\lim_{\varepsilon \to 0}  \frac{1}{2\pi i}\int_{\sigma_s}  \varphi_\varepsilon (s)  f^*(s) x^{-s} ds= \frac{1}{2\pi i}\int_{\sigma_s}  \cot(\pi s) f^*(s) x^{-s} ds$$
$$= \lim_{\varepsilon \to 0}  \frac{1}{2\pi^2 i}\int_{\sigma_s}  \left(  \int_0^{1-\varepsilon} +  \int_{1+\varepsilon}^\infty  \right) {t^{s-1}\over 1-t}\   f^*(s) x^{-s}  dt ds $$
$$= \lim_{\varepsilon \to 0}  \frac{1}{\pi} \left(  \int_0^{1-\varepsilon} +  \int_{1+\varepsilon}^\infty  \right) {f(x/t) dt \over (1-t) t } = \lim_{\varepsilon \to 0}  \frac{1}{\pi} \left(  \int_0^{x/(1+\varepsilon)} +  \int_{x/ (1-\varepsilon) }^\infty  \right) {f(t) dt \over t-x  } $$
$$= PV \   {1\over \pi}  \int_{0 }^\infty   {f(t) dt \over t-x  } = {1\over \pi} (H f)(x).$$
Therefore we proved equality (25).   Analogously we establish (26), minding that $s h(s)\in L_2(\sigma_s)$ if
 $s f^*(s) \in L_2(\sigma_s)$, where $h(s)= \cot(\pi s) f^*(s)$. Substituting these  results  in (24), it gives equality (23)
 and complete the proof of the theorem.
\end{proof}

It is well known that on $\mathbb{R}$ the following equality, involving  the iterated Hilbert transform holds, namely
$$ {1\over \pi^2} (\hat{H}^2 f)(x) =  - f(x), \quad f \in L_p(\mathbb{R}),\ p > 1,$$
where
$$  (\hat{H} f)(x)=  \int_{-\infty}^\infty \frac{f(t)}{t-x}\ dt, \quad x \in \mathbb{R}.$$
Here as an immediate corollary of equality (26) we prove the following relation between the iterated  Stieltjes and Hilbert transforms (3) and (5),
 respectively,  on $\mathbb{R}_+$, which seems to be new.  Precisely, it states

{\bf Corollary 1.} {\it Under condition $s f^*(s) \in L_2(\sigma_s)$, where $f^*$ is the Mellin
transform $(7)$ of $f$, the following relation holds for all $x >0$ }
$$ (H^2 f)(x) = (\mathcal{S}_2 f)(x) - \pi^2 f(x), \quad x > 0.$$

\begin{proof} The proof is straightforward with the use of (14), (26) and elementary trigonometric identity
$\cot^2(\pi s)= \csc^2(\pi s) -1$.

\end{proof}

{\bf Corollary 2}.  {\it Under conditions of Theorem 2 the following equality holds for convolution $(18)$}
$$ (H^2 (f*g) )(x) + \pi^2 (f*g)(x) = \sqrt x \  (\mathcal{S}_2 f)(x) (\mathcal{S}_2 g)(x), \quad x > 0.$$

\begin{proof} The proof is immediate with  the use of  factorization equality (22).

\end{proof}

Finally in this section we establish an analog of the Titchmarsh
theorem  about the absence of divisors of zero in the convolution
 (18). We have

{\bf Theorem 3}. {\it Let under conditions of Theorem $2$ $f*g= 0$.
Then either $f=0$ or $g=0$.}

\begin{proof} Indeed, one can consider the iterated  Stieltjes
transform (3) of complex variable
$$G(z)=     \int_0^\infty \frac{\log(z)- \log(t)}{z-t} f(t) dt,\quad z \in \mathbb{C} \backslash  \{0\},$$
where we take the principal branch of $\log z$. On the other hand,
we can treat $G(z)$ as the Stieltjes transform (2) of
$L_2$-function, which is analytic in  the complex plane cut along
the nonpositive real axis  (see in \cite{wid}).  In fact,
$G(z)= (\mathcal{S}_2 f)(z) = (\mathcal{S}\  (\mathcal{S} f) ) (z)$, where
similar to (11)
$$(\mathcal{S} f)(x) = {1\over 2\pi i}\int_{\sigma_s}  \Gamma(s)\Gamma(1-s) f^*(s)x^{-s} ds,  \quad x >0,$$
$$  (\mathcal{S}_2 f)(z)= {1\over 2\pi i}\int_{\sigma_s}  \left[\Gamma(s)\Gamma(1-s)\right]^2
f^*(s)z^{-s} ds, \quad |\arg z| < \pi,\  z  \neq 0$$
and  the right-hand side of the latter equality  represents  an absolutely and uniformly convergent integral in the domain
$D=\{z \in \mathbb{C},   |\arg z| < \pi,\ |z| > a >0\}.$  Indeed, with the Schwarz inequality we have 
$$\int_{\sigma_s} \left| \left[\Gamma(s)\Gamma(1-s)\right]^2
f^*(s)z^{-s} ds\right| = {\pi^2\over \sqrt{|z|}}  \int_{-\infty}^\infty  \frac{ e^{\tau\arg z} }{\cosh^2(\pi\tau)} \left| f^*\left({1\over 2}+ i\tau \right) \right| d\tau $$
$$<  {\pi^2\over \sqrt{a}} \left( \int_{-\infty}^\infty  \frac{ e^{2\pi| \tau|} }{(1/4+ \tau^2) \cosh^4(\pi\tau)} d\tau \right)^{1/2}  \left(  \int_{-\infty}^\infty   \left| \left({1\over 2}+ i\tau \right)  f^*\left({1\over 2}+ i\tau \right) \right|^2 d\tau\right)^{1/2}$$$$ < \infty.$$
Moreover,   $f(x),\   (\mathcal{S} f)(x) \in L_2(\mathbb{R}_+)$ because, evidently,  $f^*(s), \  \Gamma(s)\Gamma(1-s) f^*(s) \in
L_2(\sigma_s)$ when $s f^*(s) \in L_2(\sigma_s)$.  Thus, if $f*g=0$, then $(\mathcal{S}_2 (f*g))(z)
\equiv 0$ and by virtue of equality (22), which has a meaning for complex $z \in D$,   either $(\mathcal{S}_2
f)(z)\equiv 0$ or $(\mathcal{S}_2 g)(z)\equiv 0$ everywhere in the
complex plane cut along the nonpositive real axis.   Therefore
appealing twice to the uniqueness of the  Stieltjes transform (cf.,
for instance, in \cite{wid}, p. 336), we conclude that either $f=0$
or $g=0$ almost everywhere on $\mathbb{R}_+$.
\end{proof}

\section{A new class of singular integral equations}

This section is devoted to an application  of convolution (18) to an interesting class of integral equations, containing
a combination of the Hilbert transform and  the iterated  Stieltjes transform (or the iterated  Hilbert transform, taking into account  equality (30)).
   However, first we apply our method to simplify a solution of the singular integral equation, considered in \cite{sri},
   which involves the convolution  (4) for the Stieltjes transform (2).  Moreover, the result is known by Lemma 11.1 in \cite{sam}.
     But our  main goal will be to establish the reciprocal inverse operator, being associated with  a singular integral equation mentioned above,
      which involves   the iterated  Hilbert and Stieltjes operators.

We begin, considering  the convolution (4)  with $g(x)= x^{\alpha-1},\  0 < \alpha < 1/2$.   Namely, taking into account the value of integral (27),
we come out with the equation
$$ f(x)  \cos(\pi\alpha ) +  {\sin(\pi\alpha)\over \pi} (Hf)(x)= h(x),\  x >0,\eqno(31)$$
where  $h(x)= \pi^{-1} \sin(\pi\alpha)\  x^{1-\alpha} (f *x^{\alpha-1} )_{\mathcal{S}}$ and   $f(x),  h(x)$ satisfy conditions of Theorem 2.
 Then applying  to both sides of (31) the Mellin transform (7) and taking into account equality (25), we obtain
$$h^*(s)=  f^*(s)\left[ \cos(\pi\alpha ) +  \cot(\pi s )\sin(\pi\alpha)\right]=  f^*(s)\frac{\sin(\pi(s+\alpha))}{\sin(\pi s)}. $$
Therefore,
$$f^*(s)=  h^*(s)\frac{\sin(\pi s)} {\sin(\pi(s+\alpha))},\quad s \in \sigma_s$$
and reciprocally with the inversion formula (8) for the Mellin transform, we arrive at the unique solution of the singular integral equation (31)
$$f(x)= {1\over 2\pi i} \int_{\sigma_s} h^*(s)\frac{\sin(\pi s)} {\sin(\pi(s+\alpha))} x^{-s} ds =
\cos(\pi\alpha) h(x)-   {\sin(\pi\alpha) \over 2\pi i} $$$$\times \int_{\sigma_s} \cot(\pi(s+\alpha)) h^*(s) x^{-s} ds = \cos(\pi\alpha) h(x)-   {\sin(\pi\alpha) \over \pi} x^{\alpha}
\int_0^\infty \frac{t^{-\alpha} h(t)}{t-x }\ dt.$$
Consequently, we found a pair of reciprocal formulas for all $x >0$ and $0 < \alpha < 1/2$
$$h(x)=  \cos(\pi\alpha ) f(x) +   { \sin(\pi\alpha ) \over \pi}  \int_0^\infty \frac{f(t)}{t-x }\  dt,\eqno(32)$$
$$f(x)=   \cos(\pi\alpha)   h(x)-   {\sin(\pi\alpha)\over \pi}
\int_0^\infty \left({x\over t}\right)^\alpha \frac{ h(t)}{t-x }\  dt,\eqno(33)$$
which is  confirmed by Lemma 11.1 in \cite{sam}.

Finally in a similar manner, we apply our technique to investigate a
solvability and find a solution of a new singular integral equation,
which is associated with convolution (23) (in fact, $g(x)=
x^{\alpha-1}$ does not satisfy conditions of Theorem 2 and $(f
*x^{\alpha-1} )$ is understood by equality (23)).  Precisely,
denoting now by $h(x)= \pi^{-2} \sin^2 (\pi\alpha)\  x^{1/2 -\alpha}
(f *x^{\alpha-1} )$ and calling Corollary 1, we come out with the
following equation for all $x > 0$ and $0<\alpha < 1$
$${\cos^2(\pi\alpha)\over \pi^2}   \int_0^\infty \frac{\log(x)- \log(t)}{x-t} f(t) dt - {\sin(2\pi\alpha)\over \pi}  \int_0^\infty \frac{f(t)}{t-x }\  dt$$$$ - \cos(2\pi\alpha) f(x) = h(x).\eqno(34)$$

{\bf Theorem 4}. {\it Let $f, h$ satisfy conditions of Theorem 2. Then for all $x >0$ and $0< \alpha < 1$ singular integral equation  $(34)$ has the unique solution
$$f(x)=  {\cos^2(\pi\alpha)\over \pi^2}   \int_0^\infty \frac{\log(x)- \log(t)}{x-t}  \left({x\over t}\right)^{\alpha-1/2}  h(t) dt$$$$+  {\sin(2\pi\alpha)\over \pi}  \int_0^\infty \left({x\over t}\right)^{\alpha-1/2}  \frac{h(t)}{t-x }\  dt- \cos(2\pi\alpha) h(x).\eqno(35)$$
Conversely, singular integral equation $(35)$ has the unique solution in the form $(34)$.}

\begin{proof}   Taking the Mellin transform of both sides of (34),  minding (12), (25), (26) and Corollary 1,
after simple manipulations we  get the equality
$$h^*(s)= f^*(s) \left[ \sin(\pi\alpha) -  \cos(\pi\alpha)\cot(\pi s)\right]^2,\quad s \in \sigma_s.$$
Hence, reciprocally,
$$ f^*(s) = h^*(s) \left[\frac {\sin(\pi s)} {\sin(\pi(s+\alpha-1/2))}\right]^2=   h^*(s)
 \left[ \sin(\pi\alpha)\right.$$$$\left.  +\cos(\pi\alpha) \cot(\pi(s+\alpha-1/2))\right]^2.$$
Consequently, canceling the Mellin transform of both sides of the latter equality, we obtain
$$f(x)=  \sin^2 (\pi\alpha) h(x) +  {\sin(2\pi\alpha)\over \pi}  \int_0^\infty \left({x\over t}\right)^{\alpha-1/2}  \frac{h(t)}{t-x }\  dt$$$$+  {\cos^2 (\pi\alpha)\over \pi^2}  \  x^{\alpha-1/2} \int_0^\infty  {1\over t-x}  \int_0^\infty  \frac{u^{1/2- \alpha} \   h(u)}{u-t }\  du dt.$$
Applying again Corollary 1, we come out with solution (35).  In the  same manner we verify the  converse statement.

\end{proof}

{\bf Remark 1}. Letting $\alpha= 1/2$, we find the simplest
degenerated case of the pair of singular integral equations (34),
(35). It leads us to the unique solution $f=h$ and vice versa.

\noindent {{\bf Acknowledgments}}\\
The present investigation was supported, in part,  by the "Centro de
Matem{\'a}tica" of the University of Porto.\\

\end{document}